\documentclass[10pt,oneside]{article}
\usepackage[cp1250]{inputenc}
\usepackage[T1]{fontenc}
\usepackage{color}

\begin{document}
\newcommand{\qed}{\rule{1.5mm}{1.5mm}}
\newcommand{\proof}{\textit{Proof. }}
\newcommand{\ccon}{\rightarrowtail}
\newtheorem{theorem}{Theorem}[section]
\newtheorem{lemma}[theorem]{Lemma}
\newtheorem{remark}[theorem]{Remark}
\newtheorem{example}[theorem]{Example}
\newtheorem{corollary}[theorem]{Corollary}
\newtheorem{proposition}[theorem]{Proposition}

\begin{center}
{\LARGE\textbf{Approximation of maps into spheres by piecewise-regular maps of class $C^k$}\vspace*{3mm}}\\
{\large Marcin Bilski\footnote{M. Bilski: Department of Mathematics and Computer Science,
Jagiellonian University, \L ojasiewicza 6, 30-348 Krak\'ow,
Poland. E-mail: Marcin.Bilski@im.uj.edu.pl\vspace*{1mm}\\
Research partially supported by the NCN grant 2014/13/B/ST1/00543}}\vspace*{8mm}\\
\end{center}
\begin{abstract}\noindent The aim of this paper is to prove that every continuous map from\linebreak a compact subset of a real algebraic
variety into a sphere
can be approximated by piecewise-regular maps of class $\mathcal{C}^k,$ where $k$ is an arbitrary nonnegative integer.\vspace*{1mm}\\
\textbf{Keywords: } real algebraic variety, continuous rational map, regulous map, quasi-regulous map, approximation.\vspace*{1mm}\\
\textbf{MSC: } 14P05, 14P10, 26C15
\end{abstract}
\section{Introduction}
\label{sectionintro}
The problem of algebraic approximation of continuous maps into spheres has been studied for
many years (cf. \cite{BoKu1987}, \cite{Ku2014}, \cite{BCR} and references
therein).

Since regular maps are often to rigid to approximate arbitrary continuous maps (cf.
\cite{BoKu1987b}, \cite{BoKu1988}, \cite{BCR}), it is natural to ask whether approximation by
maps from larger classes is possible. One of such classes is the class of continuous rational
maps (in real algebraic geometry studied systematically for the first time in \cite{Ku2009}).
On nonsingular real algebraic varieties it coincides with the class of regulous maps (also
known as continuous hereditarily rational maps cf. \cite{FHMM}, \cite{KoNo}). Continuous
rational maps and regulous maps have attracted a lot of attention in recent years (see
\cite{FHMM}, \cite{FMQ}, \cite{KoKuKu}, \cite{KoNo}, \cite{Ku2014}, \cite{KuKu2},
\cite{MoJP}, \cite{Zie} and references therein). It has turned out, for example, that every
continuous map between spheres can be approximated by continuous rational maps (see
\cite{Ku2014}). However, not every continuous map from an arbitrary compact nonsingular real
algebraic variety into a sphere can be approximated by continuous rational maps (see also
\cite{Ku2014}).

In this paper (as in \cite{Ku2014}, \cite{Ku2017}), a real algebraic variety is, by
definition, a locally ringed space isomorphic to an algebraic subset of $\mathbf{R}^m,$ for
some $m,$ endowed with the Zariski topology and the sheaf of real-valued regular functions.
Each real algebraic variety carries also the Euclidean topology induced by the standard
metric in $\mathbf{R}.$ Unless explicitly stated otherwise, all topological notions relating
to real algebraic varieties refer to the Euclidean topology.

In a recent paper \cite{Ku2017}, Kucharz introduced a class of piecewise-regular maps. We say
that a continuous map $f:V\supset X\rightarrow W,$ where $V$ and $W$ are real algebraic
varieties and $X\subset V$ is some subset, is a piecewise-regular map if the following holds.
There is a stratification $\mathcal{S}$ of $V$ such that for every stratum $S\in\mathcal{S}$
the restriction of $f$ to each connected component of $S\cap X$ is a regular map. (For
precise definition see Section \ref{DWA}.)

Under the assumption that $X$ is a compact subset of a real algebraic variety, Kucharz proved
that every continuous map from $X$ into a sphere can be approximated by piecewise-regular
maps (cf. \cite{Ku2017}, Theorems 1.5). Another result of \cite{Ku2017} (see Theorem 1.6)
says that if $V$ is a compact nonsingular real algebraic variety, then each continuous map
from $V$ into a sphere is homotopic to a piecewise-regular map of class $\mathcal{C}^k,$
where $k$ is an arbitrary nonnegative integer. (For $k=\infty$ this does not hold true, see
\cite{Ku2017}, Ex. 1.7.)

The following question concerning generalization of Theorems 1.5, 1.6 of \cite{Ku2017} has
been asked in \cite{Ku2017} (see the paragraph preceding Example 1.8). Let $V$ be a compact
nonsingular real algebraic variety and $m, k,$ be positive integers. Can every continuous map
from $V$ into the $m$-dimensional unit sphere be approximated by piecewise-regular maps of
class $\mathcal{C}^k$? The aim of the present paper is to answer this question in
affirmative. To do this, we introduce a new concept of quasi-regulous map. It turns out that
the class of quasi-regulous maps is a subclass of the class of piecewise-regular maps. We
show that continuous maps from $V$ into the $m$-dimensional unit sphere can be approximated
by quasi-regulous maps of class $\mathcal{C}^k.$

{Let $W$ be a Zariski open subset of $\mathbf{R}^n.$ Recall} that a map
{$g:W\rightarrow\mathbf{R}^l$} is said to be a regulous map if it is continuous and there are
polynomials $P_1,\ldots, P_l, Q_1,\ldots, Q_l\in\mathbf{R}[T_1,\ldots, T_n]$ such that no
$Q_i$ is identically equal to {zero and}
$$g(x)=(\frac{P_1(x)}{Q_1(x)},\ldots,\frac{P_l(x)}{Q_l(x)})$$ for all points {$x\in W$} with
$Q_1(x)\neq 0,\ldots, Q_l(x)\neq 0$. ({As mentioned above, on nonsingular real algebraic
varieties "regulous" means "continuous rational". These maps have been investigated in more
general context (cf. \cite{FHMM}, \cite{KoNo}, \cite{Ku2014}), but here it is sufficient to
consider them on Zariski open subsets of $\mathbf{R}^n$.})
\vspace*{2mm}\\
{\textbf{Definition} Let $W$ be a Zariski open subset of $\mathbf{R}^n$ and let $k$ be a
nonnegative integer. A function $g:W\rightarrow\mathbf{R}$ is called $\mathcal{C}^k$
quasi-regulous if it is of class $\mathcal{C}^k$ and there is a regulous function
$f:W\rightarrow\mathbf{R}$ such that $|g(x)|=|f(x)|,$ for every $x\in W.$ Let $M$ be any (not
necessarily algebraic) subset of $\mathbf{R}^n.$\linebreak A function
$\tilde{g}:M\rightarrow\mathbf{R}$ is called $\mathcal{C}^k$ quasi-regulous if there are a
Zariski open neighborhood $W$ of $M$ in $\mathbf{R}^n$ and a $\mathcal{C}^k$ quasi-regulous
function $g:W\rightarrow \mathbf{R}$ such that $\tilde{g}={g}|_M.$ A map
$h=(h_1,\ldots,h_l):M\rightarrow \mathbf{R}^l$ is called $\mathcal{C}^k$ quasi-regulous  if
$h_i$ is a $\mathcal{C}^k$ quasi-regulous function, for
every $i=1,\ldots,l.$} \vspace*{2mm}\\
{\textit{Convention.} By a quasi-regulous map we mean a $\mathcal{C}^0$ quasi-regulous
map.}\vspace*{2mm}

Since real algebraic varieties are, by definition, isomorphic to real algebraic sets, we may
assume that $V$ in the question stated above is a real algebraic subset of some
$\mathbf{R}^n$. Actually, we prove that every continuous map from any compact (not
necessarily algebraic) subset $M$ of $\mathbf{R}^n$ to the unit sphere
$\mathbf{S}^m\subset\mathbf{R}^{m+1}$ can be approximated by $\mathcal{C}^k$ quasi-regulous
maps. As in \cite{Ku2017}, approximation is expressed in terms of the compact-open topology
of the space $\mathcal{C}(M,\mathbf{S}^m)$ of all continuous maps from $M$ to $\mathbf{S}^m$.
More precisely, by definition, a map $f\in\mathcal{C}(M,\mathbf{S}^m)$ can be approximated by
$\mathcal{C}^k$ quasi-regulous maps if every neighborhood of $f$ (with respect to the
compact-open topology) in $\mathcal{C}(M,\mathbf{S}^m)$ contains a $\mathcal{C}^k$
quasi-regulous map.

{In view of Corollary \ref{quasipiece}} (see Section \ref{DWA}), the affirmative answer to
the question follows immediately from Theorem \ref{main} which is our main result.
\label{intro}
{\begin{theorem}\label{main}Let $M\subset\mathbf{R}^n$ be any compact subset
and let $\mathbf{S}^m\subset\mathbf{R}^{m+1}$ be the $m$-dimensional unit sphere. Then for
every positive integer $k,$ every continuous map from $M$ to $\mathbf{S}^m$ can be
approximated by $\mathcal{C}^k$ quasi-regulous maps.
\end{theorem}}

The notion of $\mathcal{C}^k$ quasi-regulous map defined above on an arbitrary subset of any
$\mathbf{R}^n$ has a natural extension to maps defined on arbitrary subsets of real algebraic
varieties (via isomorphisms between real algebraic varieties and real algebraic subsets of
$\mathbf{R}^n$). Hence, in Theorem \ref{main}, the space $\mathbf{R}^n$ can be replaced by a
real algebraic variety.

The organization of this paper is as follows. In Section \ref{DWA} we present some
preliminary results on piecewise-regular maps and quasi-regulous maps and check the inclusion
between these classes. Here we also collect basic facts on differential properties of
quasi-regulous maps. In Section \ref{TRZY} the proof of Theorem \ref{main} is given.

\section{Preliminaries}\label{DWA}
The following generalization of the definition of regular map introduced in \cite{Ku2017}
will be useful to us.\vspace*{2mm}\\
\textbf{Definition} Let $V, W$ be real algebraic varieties, $X\subset V$ some (nonempty)
subset, and $Z$ the Zariski closure of $X$ in $V.$ A map $f:X\rightarrow W$ is said to be
regular if there is a Zariski open neighborhood $Z_0\subseteq Z$ of $X$ and a regular map
$\tilde{f}: Z_0\rightarrow W$ such that $\tilde{f}|_X=f.$\vspace*{2mm}

By a \textit{stratification} of a real algebraic variety $V$ we mean a finite collection of
pairwise disjoint Zariski locally closed subvarieties (some possibly empty) whose union is
equal to
$V.$\vspace*{0mm}\\
\\
\textbf{Definition} Let $V, W$ be real algebraic varieties, $f:X\rightarrow W$ a continuous
map defined on some subset $X\subset V,$ and $\mathcal{S}$ a stratification of $V.$ The map
$f$ is said to be piecewise $\mathcal{S}$-regular if for every stratum $S\in\mathcal{S}$ the
restriction of $f$ to each connected component of $X\cap S$ is a regular map (when $X\cap S$
is non-empty). Moreover, $f$ is said to be piecewise-regular if it is piecewise
$\mathcal{T}$-regular for some stratification $\mathcal{T}$ of $V.$\vspace*{2mm}

The next theorem requires the notion of \textit{nonsingular algebraic arc} (cf.
\cite{Ku2017}). Given a real algebraic variety $V,$ a subset $A\subset V$ is said to be a
nonsingular algebraic arc if its Zariski closure $C$ in $V$ is an algebraic curve (that is,
$\mathrm{dim}(C)=1$), $A\subset C\setminus\mathrm{Sing}(C),$ and $A$ is homeomorhpic to
$\mathbf{R}$.

The following result is taken from \cite{Ku2017} (Theorem 2.9).
\begin{theorem}\label{piecechar} Let $V, W$ be real algebraic varieties, $X\subset V$ a
semialgebraic subset
and $f:X\rightarrow W$ a continuous semialgebraic map. Then the following conditions are equivalent:\vspace*{2mm}\\
(a) The map $f$ is piecewise-regular\\
(b) For every nonsingular algebraic arc $A$ in $V$ with $A\subset X,$ there exists a nonempty
open subset $A_0\subset A$ such that $f|_{A_{0}}$ is a regular map.
\end{theorem}
Now we prove a corollary which explains the relation between quasi-regulous maps and
piecewise-regular maps.
\begin{corollary}\label{quasipiece}Let {$M\subset\mathbf{R}^n$ be any}
subset. Every quasi-regulous map from {$M$} to some $\mathbf{R}^m$ is a piecewise-regular
map.
\end{corollary}
\proof {Let $\tilde{f}:M\rightarrow\mathbf{R}^m$ be a quasi-regulous map. By definition,
there is a quasi-regulous map $f=(f_1,\ldots,f_m):X\rightarrow\mathbf{R}^m$ defined on a
Zariski open neighborhood $X$ of $M$ in $\mathbf{R}^n$ such that $f|_M=\tilde{f}.$ We show
that $f$ is piecewise-regular which clearly implies that $\tilde{f}$ is piecewise-regular
too.}

{We check that $f$ is regular on some open subsets of nonsingular algebraic arcs and then
apply Theorem \ref{piecechar}. Let $A\subset X$ be a nonsingular algebraic arc and $C$ be the
Zariski closure of $A$ in $X$. Let $g=(g_1,\ldots,g_m):X\rightarrow\mathbf{R}^m$ be a
regulous map such that $|f_i(x)|=|g_i(x)|$ for every $i\in\{1,\ldots,m\}$ and $x\in X.$ }

{Zariski open sets in $\mathbf{R}^n$ are nonsingular real algebraic varieties hence by
Proposition 7 of \cite{KoNo} there are $P_i,$ $Q_i\in\mathbf{R}[T_1,\ldots, T_n]$ such that
$g_i|_C=\frac{P_i}{Q_i}|_C$ outside the zero-set of $Q_i$ and $Q_i$ is not identically equal
to zero on $C,$ for $i=1,\ldots,m.$ Set $E_i:=(Q_i|_C)^{-1}(0)\cup(P_i|_C)^{-1}(0)$ and
observe that, for every $i,$ either $A\cap E_i$ is finite or $A\subset E_i$. It is clear that
in the second case, $f_i|_A=0.$ In the first one, $f_i|_{A}$ is regular on every connected
component of $A\setminus E_i.$ This is because on such a component $S$, $f_i$ has constant
sign (by continuity) so, either $f_i|_S=\frac{P_i}{Q_i}|_S$ or $f_i|_S=-\frac{P_i}{Q_i}|_S,$
and $Q_i$ vanishes nowhere on $S.$ Therefore the map $f$ is regular on every connected
component of $A\setminus T,$ where $T$ is the union of all $E_i$ such that $A$ is not a subset of $E_i.$
Now the proof is complete by Theorem~\ref{piecechar}.}\qed\vspace*{2mm}

The following lemmas will be used in the proof of our main result.

\begin{lemma}\label{tech}Let $l, k\in\mathbf{N},$ $l, k\geq 1.$\vspace*{2mm}\\
(a) Let $f(y_1,\ldots,y_p)=(y_1^{2kl}+\ldots+y_p^{2kl})^{\frac{1}{l}}.$ Then
$f\in\mathcal{C}^{k}(\mathbf{R}^p)$ and every partial derivative of $f$ of order
up to $k$ vanishes at $0\in\mathbf{R}^p$.\vspace*{2mm}\\
(b) Let $g(x,y)=\frac{x^{2k}y^{2k}}{(x^{2kl}+y^{2kl})^{\frac{1}{l}}}.$ Then
$g\in\mathcal{C}^{k}(\mathbf{R}^2)$ and every partial derivative of $g$ of order
up to $k$ vanishes at $0\in\mathbf{R}^2$.\vspace*{2mm}\\
(c) Let $h(x,y)=\frac{x^{4k}}{(x^{2kl}+y^{2kl})^{\frac{1}{l}}}.$ Then
$h\in\mathcal{C}^{k}(\mathbf{R}^2)$ and every partial derivative of $h$ of order up to $k$
vanishes at $0\in\mathbf{R}^2$.
\end{lemma}
\proof $(a)$ Fix any $j\in\mathbf{N}$ with $0\leq j\leq k.$ Observe that every partial
derivative of $f$ at $(y_1,\ldots,y_p)\neq (0,\ldots,0)$ of order $j$ is the sum of terms of
the following form, multiplied by constants:
$$(y_1^{2kl}+\ldots+y_p^{2kl})^{\frac{1}{l}-t}y_1^{(2kl-1)s_1-m_1}\cdot\ldots\cdot y_p^{(2kl-1)s_p-m_p},$$
where $t, s_i, m_i\in\mathbf{N}$ satisfy $0\leq t\leq j$ and $s_1+\ldots+s_p=t,$ and
$m_1+\ldots+m_p=j-t.$ Moreover, $(2kl-1)s_i-m_i\geq 0,$ for $i=1,\ldots,p.$ (By a partial
derivative of $f$ of order zero we mean $f$ itself.) Clearly, every such term is an analytic
function in $(y_1,\ldots,y_p)$ on $\mathbf{R}^p\setminus\{0\}^p.$ We will check that every
such term tends to zero as $(y_1,\ldots,y_p)$ approaches zero. This is an immediate
consequence of the following inequalities.

$$|(y_1^{2kl}+\ldots+y_p^{2kl})^{\frac{1}{l}-t}y_1^{(2kl-1)s_1-m_1}\cdot\ldots\cdot
y_p^{(2kl-1)s_p-m_p}|\leq$$
$$\leq(y_1^{2kl}+\ldots+y_p^{2kl})^{\frac{1}{l}-t}(\max_i|y_i|)^{(2kl-1)t-j+t}\leq$$
$$\leq p(\max_i|y_i|)^{2kl\cdot(\frac{1}{l}-t)}(\max_i|y_i|)^{(2klt-j)}=p(\max_i|y_i|)^{(2k-j)}.$$
This implies that every partial derivative of $f|_{\mathbf{R}^p\setminus\{0\}^p}$ at
$(y_1,\ldots,y_p)$ of order up to $k$ tends to zero as $(y_1,\ldots,y_p)$ approaches zero.

In order to see that every partial derivative of $f$ at $(0,\ldots,0)$ of order up to $k$
exists and equals $0$, it is sufficient to observe that when we divide the righthand side of
the last inequality, for $0\leq j\leq (k-1)$, by $\max_i|y_i|,$  the obtained expression
still tends to zero as $(y_1,\ldots,y_p)$ approaches zero. The proof of $(a)$ is
complete.\vspace*{2mm}\\
$(b)$ Fix any $j\in\mathbf{N}$ with $0\leq j\leq k.$ Observe that every partial derivative of
$g$ at $(x,y)\neq(0,0)$ of order $j$ is the sum of terms of the following form, multiplied by
constants:
$$(x^{2kl}+y^{2kl})^{-\frac{1}{l}-t}x^{(2kl-1)s_1+2k-m_1}\cdot y^{(2kl-1)s_2+2k-m_2},$$
where $t, s_i, m_i\in\mathbf{N}$ satisfy $0\leq t\leq j$ and $s_1+s_2=t,$ and $m_1+m_2=j-t.$
Moreover, $(2kl-1)s_i+2k-m_i\geq 0,$ for $i=1,2.$ Clearly, every such term is an analytic
function in $(x,y)$ on $\mathbf{R}^2\setminus\{0\}^2.$ We will check that every such term
tends to zero as $(x,y)$ approaches zero. This is an immediate consequence of the following
inequalities.
$$|(x^{2kl}+y^{2kl})^{-\frac{1}{l}-t}x^{(2kl-1)s_1+2k-m_1}\cdot y^{(2kl-1)s_2+2k-m_2}|\leq$$
$$\leq(x^{2kl}+y^{2kl})^{-\frac{1}{l}-t}\cdot\max\{|x|,|y|\}^{(2kl-1)t+4k-j+t}\leq$$
$$\max\{|x|,|y|\}^{2kl\cdot(-\frac{1}{l}-t)}\cdot\max\{|x|,|y|\}^{(2klt+4k-j)}=
\max\{|x|,|y|\}^{2k-j}.$$ This implies that every partial derivative of
$g|_{\mathbf{R}^2\setminus\{0\}^2}$ at $(x,y)$ of order up to $k$ tends to zero as $(x,y)$
approaches zero.

In order to see that every partial derivative of $g$ at $(0,0)$ of order up to $k$ exists and
equals $0,$ it is sufficient to observe that when we divide the righthand side of the last
inequality, for $0\leq j\leq (k-1)$, by $\max\{|x|,|y|\},$  the obtained expression still
tends to
zero as $(x,y)$ approaches zero. The proof of $(b)$ is complete.\vspace*{2mm}\\
$(c)$ can be proved in the way very similar to $(b)$. \qed\vspace*{2mm}

The class of functions defined below will be used for several times in the
sequel.\vspace*{2mm}\\
\textbf{Definition} For any {open subset $U$} of $\mathbf{R}^n,$ let $\mathcal{C}^k_{l}(U)$
denote the class of all functions $v:U\rightarrow\mathbf{R}$ for which
$|v|^{\frac{1}{l}}\in\mathcal{C}^k(U).$

\begin{lemma}\label{descent}Let {$U$ be an open subset} of $\mathbf{R}^n.$
Let $v_1,\ldots,v_p\in\mathcal{C}^k_{2kl}(U)$ be nonnegative functions, where
$k,l\in\mathbf{N}$ and $k,l\geq 1$. Then $v_1+\ldots+v_p\in\mathcal{C}^k_l(U)$ and
$\frac{v_1v_2}{v_1+v_2}\in\mathcal{C}^k_l(U),$ and
$\frac{v_1^2}{v_1+v_2}\in\mathcal{C}^k_l(U).$ Here, if $v_1(x)=v_2(x)=0,$ then, by
definition, $\frac{v_1(x)v_2(x)}{v_1(x)+v_2(x)}=\frac{v_1^2(x)}{v_1(x)+v_2(x)}=0.$
\end{lemma}
\proof There are $u_1,\ldots,u_p\in\mathcal{C}^k(U)$ such that
$v_1=u_1^{2kl},\ldots,v_p=u_p^{2kl}.$

The map $(v_1+\ldots+v_p)^{\frac{1}{l}}$ is the composition of the map
$f(y_1,\ldots,y_p)=(y_1^{2kl}+\ldots+y_p^{2kl})^{\frac{1}{l}}$ with
$x\mapsto(u_1(x),\ldots,u_p(x)).$ By Lemma \ref{tech} $(a)$, we have
$(v_1+\ldots+v_p)^{\frac{1}{l}}\in\mathcal{C}^k(U),$ hence
$v_1+\ldots+v_p\in\mathcal{C}^k_l(U)$.

Similarly, the map $(\frac{v_1v_2}{v_1+v_2})^{\frac{1}{l}}$ is the composition of
$g(y_1,y_2)=\frac{y_1^{2k}y_2^{2k}}{(y_1^{2kl}+y_2^{2kl})^{\frac{1}{l}}}$ with
$x\mapsto(u_1(x),u_2(x)),$ and the map $(\frac{v_1^2}{v_1+v_2})^{\frac{1}{l}}$ is the
composition of $h(y_1,y_2)=\frac{y_1^{4k}}{(y_1^{2kl}+y_2^{2kl})^{\frac{1}{l}}}$ with
$x\mapsto(u_1(x),u_2(x)).$ Therefore, by Lemma \ref{tech} $(b)$, $(c)$, the proof is
complete.\qed
\begin{lemma}\label{symmetr}Let {$U$ be an open subset} of $\mathbf{R}^n.$
Let $f\in\mathcal{C}^{k}_{l}(U)$, where $l, k\in\mathbf{N}$ with $k\geq 1$ and $l\geq k+1.$
Then for every continuous function $g:U\rightarrow\mathbf{R}$ such that for every $x\in U$
either $g(x)=f(x)$ or $g(x)=-f(x)$, we have $g\in\mathcal{C}^k(U).$
\end{lemma}
\proof By hypothesis, {there is $h\in\mathcal{C}^k(U)$} such that $h^l=|f|=|g|.$ In
particular, $h^{-1}(0)=f^{-1}(0)=g^{-1}(0).$ By continuity of $g,$ for every connected
component $A_i$ of $U\setminus(h^{-1}(0))$ there is $\beta_i\in\{-1,1\}$ with $g=\beta_i\cdot
h^l$ on $A_i$ and $g|_{h^{-1}(0)}=0.$

Now, it is not difficult to observe that for such $g$ there exist all partial derivatives of
order up to $k$ at every point of $U.$ These partial derivatives vanish identically on
$h^{-1}(0)$ and are continuous on $U.$\qed\vspace*{2mm}

Let $N$ be any {compact subset of $\mathbf{R}^n$ and let $h:N\rightarrow\mathbf{R}$} be a
continuous function. For any $\delta>0$ put {$N_{\delta}=\{x\in\mathbf{R}^n: dist(x,
N)<\delta\}$} and put $||h||_N=\mathrm{sup}_{x\in N}|h(x)|.$ We assume
$\emptyset_{\delta}=\emptyset$ and $||h||_{\emptyset}=0.$

\begin{lemma}\label{level} Let {$G\subset\subset\mathbf{R}^n$ be an open
ball} and let $f,g:\overline{G}\rightarrow\mathbf{R}$ be continuous, nonnegative,
semialgebraic functions. {Assume that $(f\cdot g)^{-1}(0)$ is a compact
$\mathcal{C}^{\infty}$ submanifold of $G$ of pure codimension $1$.} Then for every
$\varepsilon>0$ there are $P,Q\in\mathbf{R}[x_1,\ldots,x_n]$ and an open tubular neighborhood
$T\subset\subset G$ of $(f\cdot g)^{-1}(0)$ of constant radius such that $||f\cdot
g||_{\overline{T}}<\varepsilon$, {$(P\cdot Q)^{-1}(0)\cap\partial G=\emptyset$}
and such that the following conditions are satisfied:\vspace*{2mm}\\
(a) for every connected component $A$ of $G\setminus (f\cdot g)^{-1}(0)$ there is a connected
component $B$ of $G\setminus (P\cdot Q)^{-1}(0)$ with $A\subset B\cup\overline{T}.$ Moreover,
for every connected component $B$ of $G\setminus (P\cdot Q)^{-1}(0)$, there is at most one
connected component
$A$ of $G\setminus (f\cdot g)^{-1}(0)$ with $A\subset B\cup\overline{T}.$\vspace*{2mm}\\
(b) for every $l\in\mathbf{N}$ there are $\tilde{f}, \tilde{g}\in\mathbf{R}[x_1,\ldots,x_n]$,
$\tilde{f}, \tilde{g}> 0$ on $\mathbf{R}^n,$ such that $||f-\tilde{f}\cdot P^{2l}||_G\leq
\varepsilon, ||g-\tilde{g}\cdot Q^{2l}||_G\leq \varepsilon.$
\end{lemma}
\textit{Remark} The algebraicity assumptions above can be relaxed, but we prefer to state
this lemma in the
setting in which it will be used.\vspace*{2mm}\\
\textit{Proof of Lemma \ref{level} }Fix $\varepsilon >0$ (we may assume $\varepsilon<1$).
Choose an open tubular neighborhood $T$ of $(f\cdot g)^{-1}(0)$ {relatively compact} in $G$
of constant radius so small that $||f\cdot g||_{\overline{T}}<\varepsilon$ and for every
connected component {$A$ of $G\setminus (f\cdot g)^{-1}(0)$} there is precisely one connected
component $\hat{A}$ of $G\setminus\overline{T}$ such that $\hat{A}\subset A$.

Next choose $\delta>0$ so small that $(f^{-1}(0))_{2\delta}\cup(g^{-1}(0))_{2\delta}\subset
T,$ and moreover, $||f||_{(f^{-1}(0))_{2\delta}}<\frac{\varepsilon}{2}$ and
$||g||_{(g^{-1}(0))_{2\delta}}<\frac{\varepsilon}{2}$.

Now there is $0<\varepsilon_1<1$ such that $\{x\in\overline{G}
:f(x)\leq{\varepsilon_1}\}\subset(f^{-1}(0))_{\delta}$ and $\{x\in\overline{G}
:g(x)\leq{\varepsilon_1}\}\subset(g^{-1}(0))_{\delta}.$ Using the Stone-Weierstrass theorem,
approximate $f-\frac{{\varepsilon_1}}{2}$ and $g-\frac{{\varepsilon_1}}{2}$ on
{$\overline{G}$} by polynomials $P, Q,$ respectively, with precision high enough to ensure
that {$P^{-1}(0)\cap\overline{G}\subset(f^{-1}(0))_{\delta}$,
$Q^{-1}(0)\cap\overline{G}\subset(g^{-1}(0))_{\delta},$ $P|_{f^{-1}(0)}<0,$
$Q|_{g^{-1}(0)}<0,$ $P|_{\overline{G}\setminus (f^{-1}(0))_{\delta}}>0,$
$Q|_{\overline{G}\setminus (g^{-1}(0))_{\delta}}>0.$} In particular, it follows that
{$(P\cdot Q)^{-1}(0)\cap G\subset T$} and {$(P\cdot Q)^{-1}(0)\cap\partial G=\emptyset$}.

{Proof of $(a)$.} Let $\hat{A}_1,\ldots,\hat{A}_s$ denote all pairwise distinct connected
components of {$G\setminus\overline{T}.$} By the previous paragraph, $\hat{A}_j\cap (P\cdot
Q)^{-1}(0)=\emptyset,$ so there is the uniquely determined connected component $B_j$ of
{$G\setminus(P\cdot Q)^{-1}(0)$} with $\hat{A}_j\subset B_j,$ for $j=1,\ldots,s.$

We shall need the following property of $\{B_i\}.$  If $i\neq j,$ then $B_i\neq B_j.$ If this
was not true, then there were $a_i\in \hat{A}_i$ and $a_j\in \hat{A}_j$ and a path connecting
these points in $B_i=B_j.$ The path must intersect $(f\cdot g)^{-1}(0)$ at some point $b$
because $\hat{A}_i, \hat{A}_j$ are contained in two distinct connected components of
{$G\setminus(f\cdot g)^{-1}(0)$}. Then either $f(b)=0$ or $g(b)=0$. Assume that $f(b)=0.$
(For $g(b)=0$ we proceed identically.) We have $P(b)<0$ and $P(a_i)>0.$ Consequently, at some
point $d$ of the path $P(d)=0.$ This implies that $B_i\cap (P\cdot Q)^{-1}(0)\neq\emptyset,$
a contradiction.

Now we show that for every connected component $B$ of {$G\setminus (P\cdot Q)^{-1}(0)$},
there is at most one connected component $A$ of {$G\setminus (f\cdot g)^{-1}(0)$} such that
$A\subset B\cup\overline{T}.$ Fix $B.$ Observe that for every $i\in\{1,\ldots,s\}$ either
$\hat{A}_i\subset B$ or $\hat{A}_i\cap B=\emptyset.$ Indeed, if for some $i$ we had
$\hat{A}_i\cap B\neq\emptyset$ and $\hat{A}_i$ was not a subset of $B,$ then from the connectedness of
$\hat{A}_i$ we had $\partial B\cap\hat{A}_i\neq\emptyset.$ But this is impossible because
$(P\cdot Q)|_{\partial B\cap G}=0$ and $(P\cdot Q)^{-1}(0)\cap \hat{A}_i=\emptyset.$

Suppose that $\hat{A}_i\cup\hat{A}_j\subset B\cup\overline{T}$ for $i\neq j.$ We know that
$\hat{A}_i\cup\hat{A}_j$ is not a subset of $B$ so, in view of previous paragraph, either
$\hat{A}_i\subset\overline{T}$ or $\hat{A}_j\subset\overline{T},$ a contradiction.

Let us verify the first claim of $(a)$. Fix a connected component $A$ of {$G\setminus (f\cdot
g)^{-1}(0).$} Let $\hat{A}$ be the unique connected component of {$G\setminus\overline{T}$}
such that $\hat{A}\subset A.$ Since $(P\cdot Q)^{-1}(0)\cap\hat{A}=\emptyset,$ there is the
uniquely determined connected component $B$ of {$G\setminus (P\cdot Q)^{-1}(0)$} containing
$\hat{A}.$ Therefore $\hat{A}\cup \overline{T}\subset B\cup \overline{T}$ and the proof of
$(a)$ is complete because $A\subset\hat{A}\cup\overline{T}.$


Proof of $(b)$. Fix any $l\in\mathbf{N}$. Since
$||f||_{(f^{-1}(0))_{\delta}}<\frac{\varepsilon}{2}$ and
$||g||_{(g^{-1}(0))_{\delta}}<\frac{\varepsilon}{2},$ we can take $\alpha\in\mathbf{R},
\alpha>0$ so small that
$$||\frac{f+\alpha}{P^{2l}+\alpha}\cdot
P^{2l}||_{(f^{-1}(0))_{\delta}}<\frac{\varepsilon}{2} \mbox{ and }
||\frac{g+\alpha}{Q^{2l}+\alpha}\cdot
Q^{2l}||_{(g^{-1}(0))_{\delta}}<\frac{\varepsilon}{2}.$$ Since {$P^{-1}(0)\cap G
\subset(f^{-1}(0))_{\delta}$} and {$Q^{-1}(0)\cap G\subset(g^{-1}(0))_{\delta}$} we can
additionally require that {$$||\frac{f+\alpha}{P^{2l}+\alpha}\cdot
P^{2l}-f||_{G\setminus(f^{-1}(0))_{\delta}}<\frac{\varepsilon}{2} \mbox{ and }
||\frac{g+\alpha}{Q^{2l}+\alpha}\cdot
Q^{2l}-g||_{G\setminus(g^{-1}(0))_{\delta}}<\frac{\varepsilon}{2}.$$} Note that $f+\alpha,
g+\alpha>0,$ hence for close polynomial approximations $\bar{f}, \bar{g}$ of
{$\sqrt{\frac{f+\alpha}{P^{2l}+\alpha}}, \sqrt{\frac{g+\alpha}{Q^{2l}+\alpha}}$ on
$\overline{G},$ after defining $\tilde{f}:=\bar{f}^2+\eta, \tilde{g}:=\bar{g}^2+\eta$, for
small $\eta\in\mathbf{R}, \eta>0,$ we have polynomial approximations $\tilde{f}, \tilde{g}$
of $\frac{f+\alpha}{P^{2l}+\alpha},$ $\frac{g+\alpha}{Q^{2l}+\alpha}$ on $\overline{G}$ with
$\tilde{f}>0,\tilde{g}>0$ on $\mathbf{R}^n$ and}
$$||\tilde{f}\cdot P^{2l}||_{(f^{-1}(0))_{\delta}}<\frac{\varepsilon}{2}, \mbox{ }
||\tilde{g}\cdot Q^{2l}||_{(g^{-1}(0))_{\delta}}<\frac{\varepsilon}{2},$$ and
{$$||\tilde{f}\cdot P^{2l}-f||_{G\setminus(f^{-1}(0))_{\delta}}<\frac{\varepsilon}{2},\mbox{
} ||\tilde{g}\cdot Q^{2l}-g||_{G\setminus(g^{-1}(0))_{\delta}}<\frac{\varepsilon}{2}.$$}
Finally observe that\vspace*{3mm}\\
{$||f-\tilde{f}\cdot P^{2l}||_{G}\leq\max\{||f-\tilde{f}\cdot
P^{2l}||_{(f^{-1}(0))_{\delta}}, ||f-\tilde{f}\cdot P^{2l}||_{G\setminus
(f^{-1}(0))_{\delta}}\}\leq$}\\
$\leq\max\{||f||_{(f^{-1}(0))_{\delta}}+||\tilde{f}\cdot P^{2l}||_{(f^{-1}(0))_{\delta}},
\frac{\varepsilon}{2}\}\leq\varepsilon.$\vspace*{3mm}\\
The same argument shows that {$||g-\tilde{g}\cdot Q||_G<\varepsilon,$} proving $(b)$.\qed
\section{Proof of Theorem \ref{main}}
\label{TRZY} Let $u=(u_1,\ldots,u_{m+1}):M\rightarrow\mathbf{S}^m\subset\mathbf{R}^{m+1}$ be
a continuous map ($m\geq 1$). In particular, $u_1^2+\ldots+u_{m+1}^2=1.$ {We will show that
for every positive integer $k,$ $u$ can be approximated by $\mathcal{C}^k$ quasi-regulous
maps from $M$ to $\mathbf{S}^m\subset\mathbf{R}^{m+1}.$}

{First we show that there is a Nash (i.e. both semi-algebraic and analytic) map
$\tilde{u}=(\tilde{u}_1,\ldots,\tilde{u}_{m+1})$ from $\mathbf{R}^n\setminus\Sigma\supset M$
to $\mathbf{S}^m,$ where $\Sigma$ is a compact algebraic subset of $\mathbf{R}^n$ of
codimension at least $2,$ with the following properties. The map $\tilde{u}$ approximates $u$
on $M$ and, for every $i=1,\ldots,m+1,$ $\tilde{u}_i^{-1}(0)\cup\Sigma$ is a compact
$\mathcal{C}^{\infty}$ submanifold of $\mathbf{R}^n$ of pure dimension $n-1$ or
$\tilde{u}_i^{-1}(0)\cap M=\emptyset.$}

{To do this, apply the Stone-Weierstrass theorem to approximate $u_1,\ldots,$ $u_{m+1}$ by
polynomials $\hat{u}_1,\ldots,\hat{u}_{m+1}$ on $M$. Next choose $\varepsilon\in(0,+\infty)$
small enough and a positive integer $\hat{m}$ large enough to ensure that after having
replaced every $\hat{u}_i$ by $\hat{u}_i+\varepsilon\cdot(x_1^2+\ldots+x_n^2)^{\hat{m}}$ we
obtain a map to $\mathbf{R}^{m+1}$ with the following property. There is an open ball
$G\subset\subset\mathbf{R}^n$ such that for every $i,$ $\hat{u}_i$ is bounded from below by
$1$ outside $G$ and $\hat{u}_i$ still approximates $u_i$ on $M$.}

{Next, add small constants to all $\hat{u}_i$'s to ensure that the zero-set of every
$\hat{u}_i$ is a compact $\mathcal{C}^{\infty}$ submanifold of $G$ of pure dimension $n-1$ or
is empty and that $\mathrm{dim}(\hat{u}_1^{-1}(0)\cap\ldots\cap\hat{u}_{m+1}^{-1}(0))\leq
n-2$. Define $\Sigma:=\hat{u}_1^{-1}(0)\cap\dots\cap\hat{u}_{m+1}^{-1}(0)$ and observe that
if the approximation is close enough, then $\Sigma\cap M=\emptyset.$}

{Finally, set $\hat{u}=(\hat{u}_1,\ldots,\hat{u}_{m+1})$ and compose
$\hat{u}|_{\mathbf{R}^n\setminus\Sigma}$ with the standard retraction of
$\mathbf{R}^{m+1}\setminus\{0\}^{m+1}$ onto $\mathbf{S}^m,$ to obtain}
$$\tilde{u}=\Big{( }\frac{\hat{u}_1}{\sqrt{\hat{u}_1^2+\ldots+\hat{u}_{m+1}^2}},\ldots,
\frac{\hat{u}_{m+1}}{\sqrt{\hat{u}_1^2+\ldots+\hat{u}_{m+1}^2}}\Big{) }.$$

{To complete the proof of Theorem \ref{main} it is sufficient to approximate $\tilde{u}$ on
{$M$} by $\mathcal{C}^k$ quasi-regulous maps into $\mathbf{S}^m.$ {Hence, we may assume that
$u$ is a Nash map from $\mathbf{R}^n\setminus\Sigma\supset M$ to $\mathbf{S}^m,$ where
$\Sigma$ is a compact algebraic subset of $\mathbf{R}^n$ with $\mathrm{dim}(\Sigma)\leq
n-2.$} {Moreover, for every $i=1,\ldots,m+1,$ $u_{i}^{-1}(0)\cup\Sigma$ is a compact
$\mathcal{C}^{\infty}$ submanifold of $\mathbf{R}^n$ of pure dimension $n-1$ or
$u_i^{-1}(0)\cap M=\emptyset$.}

{If, for some $i,$ ${u}_i^{-1}(0)\cap M=\emptyset$, then ${u}|_M$ is a non-surjective map to
$\mathbf{S}^m$ which can be easily approximated by regular maps to $\mathbf{S}^m$ because the
unit sphere with a point removed is biregularly isomorphic to $\mathbf{R}^m.$ If, for every
$i,$ ${u}_i^{-1}(0)\cup\Sigma$ is a compact $\mathcal{C}^{\infty}$ submanifold of
$\mathbf{R}^n$ of pure dimension $n-1,$ then we proceed as follows.}

Let us define a new family of continuous, semialgebraic non-negative functions from
{$\mathbf{R}^n\setminus\Sigma$ to $\mathbf{R}$}:\vspace*{2mm}\\
$v_1=|u_1|,$\\
$v_{i+1}=\sqrt{v_i^2+u_{i+1}^2},\mbox{ }i=1,\ldots,m.$\vspace*{2mm}\\
In particular, $v_{m+1}=1$. It is clear that $u_i, v_i$ satisfy the following system of
equations:\vspace*{2mm}\\
$u_1^2+u_2^2=v_2^2,$\\
$v_i^2+u_{i+1}^2=v_{i+1}^2,\mbox{ for } i=2,\ldots,m.$\vspace*{2mm}\\
We will approximate $u_1$ on $M$ by polynomials and, for $i\geq 2,$ $u_i, v_i$ {on $M$} by
{$\mathcal{C}^k$ quasi-regulous functions $\bar{u}_i, \bar{v}_i$ from $\mathbf{R}^n$ to
$\mathbf{R}$} satisfying
\vspace*{2mm}\\
$\bar{u}_1^2+\bar{u}_2^2=\bar{v}_2^2,$\\
$\bar{v}_i^2+\bar{u}_{i+1}^2=\bar{v}_{i+1}^2,\mbox{ for } i=2,\ldots,m,$\vspace*{2mm}\\
{and such that $\bar{v}_{m+1}^{-1}(0)$ is an algebraic subset of $\mathbf{R}^n.$ Since
$v_{m+1}=1,$ we will have $M\cap \bar{v}^{-1}_{m+1}(0)=\emptyset$ for $\bar{v}_{m+1}$ close
to $v_{m+1}$ on $M$.} {Then, the map $w$ whose components $w_i,$ for $i=1,\ldots,m+1,$ are
defined by $w_i=\frac{\bar{u}_i}{\bar{v}_{m+1}},$} {is a $\mathcal{C}^k$ quasi-regulous map
from $\mathbf{R}^n\setminus\bar{v}_{m+1}^{-1}(0)$ to $\mathbf{S}^{m}$ approximating $u$ on
$M,$} and the proof will be completed.

Now we will prove that $u_1$ can be approximated by polynomials and, for $i\geq 2,$ $v_i$ can
be approximated by regulous functions of class $\mathcal{C}^k$ and $u_i$ can be approximated
by $\mathcal{C}^k$ quasi-regulous functions
 satisfying the equations listed above. This will be done indirectly:
approximating every $u_i,$ for $i=2,\ldots,m+1,$ will be preceded by approximating the
absolute value $|u_i|$ {on $M$} by a nonnegative regulous function
{$|\bar{u}_i|\in\mathcal{C}^{k+1}_{l}(\mathbf{R}^n),$} for some $l\geq k+1.$ (For the
definition of $\mathcal{C}^{k+1}_l(\mathbf{R}^n)$ see Section \ref{DWA}.) Using $|\bar{u}_i|$
and Lemma \ref{symmetr}, we will produce a quasi-regulous approximation of $u_i$ satisfying
the requirements.

The whole construction will be carried out in $m$ steps. In the first step we handle the
functions appearing in the equation $u_1^2+u_2^2=v_2^2$. In step number $i,$ for $i\geq 2,$
we deal with the functions appearing in $v_i^2+u_{i+1}^2=v_{i+1}^2$ making use of the data
obtained previously.

Let us discuss the first step. From the first equation we have
$|u_2|=\sqrt{v_2-u_1}\cdot\sqrt{v_2+u_1}$ {on $\mathbf{R}^n\setminus\Sigma$.} {Let
$f,g:\mathbf{R}^n\rightarrow\mathbf{R}$ be continuous nonnegative semialgebraic functions
with the following properties: $$f|_M=(\sqrt{v_2-u_1})|_M,\mbox{ }
f^{-1}(0)=(\sqrt{v_2-u_1})^{-1}(0)\cup\Sigma,$$  $$g|_M=(\sqrt{v_2+u_1})|_M,\mbox{ }
g^{-1}(0)=(\sqrt{v_2+u_1})^{-1}(0)\cup\Sigma.$$}Such functions can be obtained as follows.
By the semialgebraic version of the Urysohn lemma there are a neighborhood $U$ of $\Sigma$ in
$\mathbf{R}^n$ disjoint from $M$ and a continuous semialgebraic function $\beta:
\mathbf{R}^n\rightarrow [0,1]$ with $\beta^{-1}(0)=\Sigma$ and $\beta|_{\mathbf{R}^n\setminus
U}=1.$ Then, by Proposition 2.6.4 of \cite{BCR}, there exists an integer $N$ such that
$f=\sqrt{v_2-u_1}\cdot\beta^N,$ $g=\sqrt{v_2+u_1}\cdot\beta^N$ extended by zero on $\Sigma$
are continuous on $\mathbf{R}^n.$

{Since $(f\cdot g)^{-1}(0)=u_2^{-1}(0)\cup\Sigma$ is a compact $\mathcal{C}^{\infty}$
submanifold of $\mathbf{R}^n$ of pure dimension $n-1$, there is an open ball
$G\subset\subset\mathbf{R}^n$ such that $(f\cdot g)^{-1}(0)\cup M\subset G.$ Then, clearly,
$f|_{\overline{G}},g|_{\overline{G}}$ satisfy the hypotheses of Lemma \ref{level}.}

By Lemma \ref{level}, we obtain polynomial approximations $\hat{f}$ of $f$ and $\hat{g}$ of
$g$ {on $\overline{G}$} of the form $\hat{f}=\tilde{f}\cdot P^{2l}$ and
$\hat{g}=\tilde{g}\cdot Q^{2l}$ with any given precision $\varepsilon,$ where $P, Q,
\tilde{f}, \tilde{g}\in\mathbf{R}[x_1,\ldots,x_n],$ $\tilde{f}>0, \tilde{g}>0$ on
{$\mathbf{R}^n$} and $l=4^{m-1}(k+1)^{2m+1}.$ Define $\bar{u}_1, \bar{v}_2$ by the equations
$\hat{f}^2=\bar{v}_2-\bar{u}_1,$ $\hat{g}^2=\bar{v}_2+\bar{u}_1$ to obtain
$\bar{v}_2=\frac{\hat{f}^2+\hat{g}^2}{2},$ $\bar{u}_1=\frac{-\hat{f}^2+\hat{g}^2}{2}$ and set
$|\bar{u}_2|=\hat{f}\cdot\hat{g}.$ Then $\bar{u}_1^2+|\bar{u}_2|^2=\bar{v}_2^2.$ {Note that
$\bar{u}_1^{-1}(0)$ and $|\bar{u}_2|^{-1}(0)$ are algebraic subsets of $\mathbf{R}^n$.} It is
clear that $\bar{u}_1, |\bar{u}_2|, \bar{v}_2$ approximate $u_1, |u_2|, v_2,$ respectively
{on $M$}. To complete the first step it remains to recover $\bar{u}_2$ with certain
regularity properties when $|\bar{u}_2|$ is given. Actually, we need to control some
regularity properties of $\bar{v}_2$ as well because $\bar{v}_2$ will be used to define
$\bar{u}_3$ in the next step of the construction.

We have {$\hat{f}^2, \hat{g}^2\in\mathcal{C}^{k+1}_{2\cdot
4^{m-1}(k+1)^{2m+1}}(\mathbf{R}^n)$} which, in view of Lemma \ref{descent}, implies
{$\bar{v}_2\in\mathcal{C}^{k+1}_{4^{m-1}(k+1)^{2m}}(\mathbf{R}^n)$} and
$|\bar{u}_2|\in\mathcal{C}^{k+1}_{4^{m-1}(k+1)^{2m+1}}(\mathbf{R}^n)
\subset\mathcal{C}^{k+1}_{4^{m-1}(k+1)^{2m}}(\mathbf{R}^n).$

In order to define $\bar{u}_2,$ given $|\bar{u}_2|,$ it is sufficient to define the sign of
$\bar{u}_2$ on every connected component of {$\mathbf{R}^n\setminus(P\cdot Q)^{-1}(0).$} Let
$A_1,\ldots,A_s$ denote all pairwise distinct connected components of {$G\setminus(f\cdot
g)^{-1}(0).$} By Lemma \ref{level}, there is a tubular neighborhood {$T\subset\subset G$} of
$(f\cdot g)^{-1}(0)$ and there are connected components $B_1,\ldots,B_s$ of
{$G\setminus(P\cdot Q)^{-1}(0)$} with $A_i\subset B_i\cup\overline{T}$ and such that $B_i\neq
B_j$ if $i\neq j.$

{Recall that $u_2^{-1}(0)\cup\Sigma=(f\cdot g)^{-1}(0)$ and} set
$\mathrm{sgn}(\bar{u}_2|_{B_i})=\mathrm{sgn}(u_2|_{A_i}),$ for $i=1,\ldots,s.$ On every
component of {$G\setminus(P\cdot Q)^{-1}(0)$} different from $B_1,\ldots, B_s$ set the sign
of $\bar{u}_2$ arbitrarily and set $\bar{u}_2|_{(P\cdot Q)^{-1}(0)}=0$. {Since $(P\cdot
Q)^{-1}(0)\cap\partial G=\emptyset,$ there is precisely one connected component $C$ of
$G\setminus (P\cdot Q)^{-1}(0)$ with $\overline{C}\cap\partial G\neq\emptyset.$ Therefore on
every connected component $E$ of $(\mathbf{R}^n\setminus G)\setminus (P\cdot Q)^{-1}(0)$ we
can define $\mathrm{sgn}(\bar{u}_2|_{E})=\mathrm{sgn}(\bar{u}_2|_C$).} Clearly, $\bar{u}_2$
is a quasi-regulous function and, by Lemma \ref{symmetr},
{$\bar{u}_{2}\in\mathcal{C}^{k+1}(\mathbf{R}^n).$}

Let us check that $\bar{u}_2$ approximates $u_2$ on {$B_i\cap M,$} for $i=1,\ldots,s.$ It is
clear that $\bar{u}_2$ approximates $u_2$ on {$A_i\cap B_i\cap M.$} This is because $u_2$ and
$\bar{u}_2$ have the same sign at every point of ${A_i\cap B_i}$ and $|\bar{u}_2|$
approximates $|u_2|$ {on $M$.} On the other hand, $B_i\setminus A_i\subset\overline{T}.$ If
this was not true, then we had $((B_i\setminus A_i)\setminus\overline{T})\cap
A_j\neq\emptyset$ for some $j\neq i$ because {$\bigcup_{i=1}^sA_i\cup\overline{T}=G.$} Then
$((B_i\setminus A_i)\setminus\overline{T})\cap (B_j\cup\overline{T})\neq\emptyset$ so
$B_i\cap B_j\neq\emptyset$ for some $j\neq i,$ a contradiction. Since $B_i\setminus
A_i\subset\overline{T},$ the function $\bar{u}_2$ approximates $u_2$ on {$(B_i\setminus
A_i)\cap M$} because {$|u_2|_{M\cap\overline{T}}|=(f\cdot g)|_{M\cap\overline{T}},$ and
$|\bar{u}_2|_{M\cap\overline{T}}|$ approximates $|u_2|_{M\cap\overline{T}}|,$ and, by Lemma
\ref{level}, $(f\cdot g)|_{M\cap\overline{T}}$ may be assumed as close to zero as we wish.}
Hence, {$u_2|_{M\cap\overline{T}}$ and $\bar{u}_2|_{M\cap\overline{T}}$ are close to zero
therefore close to each other (even if $u_2$ and $\bar{u}_2$ have different signs at some
points of $M\cap\overline{T}$)}.

Let us check that $\bar{u}_2$ approximates $u_2$ on {$M\setminus\bigcup_{i=1}^s{B_i}.$}
Observe that {$G\setminus\bigcup_{i=1}^s{B_i}\subset\overline{T}.$} This is because
{$\bigcup_{i=1}^s{A_i}\cup\overline{T}={G}$} and, for every $i,$ we have ${A_i}\subset
{B_i}\cup \overline{T}.$ Consequently, $M\setminus\bigcup_{i=1}^s{B_i}\subset\overline{T}$
and the assertion follows as above. Thus we have proved that $\bar{u}_2$ approximates $u_2$
on {$M.$}

Steps numbered by $i\in\{2,\ldots, m\}$ are slightly different. We will describe them
inductively. Assume that step number $j$, where $1\leq j <m,$ has been accomplished i.e. we
have continuous functions $\bar{u}_{j+1}, \bar{v}_{j+1}$ approximating $u_{j+1},$ $v_{j+1},$
respectively, {on $M$} such that $|\bar{u}_{j+1}|,$
{$\bar{v}_{j+1}\in\mathcal{C}^{k+1}_{4^{m-j}(k+1)^{2(m-j+1)}}(\mathbf{R}^n)$} are
non-negative regulous functions satisfying $\bar{v}_{j}^2+\bar{u}_{j+1}^2=\bar{v}_{j+1}^2$
(we assume {$\bar{v}_1=|\bar{u}_1|$}). In particular, $\bar{u}_{j+1}$ is a quasi-regulous
function and, by Lemma \ref{symmetr}, {$\bar{u}_{j+1},
\bar{v}_{j+1}\in\mathcal{C}^{k+1}(\mathbf{R}^n).$}

We will define continuous functions $\bar{u}_{j+2},$ $\bar{v}_{j+2}$ approximating
$u_{j+2}, v_{j+2}$, respectively, {on $M$}
 such that $|\bar{u}_{j+2}|,$ $\bar{v}_{j+2}\in$
$\mathcal{C}^{k+1}_{4^{m-j-1}(k+1)^{2(m-j)}}(\mathbf{R}^n)$ are non-negative regulous
functions satisfying $\bar{v}_{j+1}^2+\bar{u}_{j+2}^2=\bar{v}_{j+2}^2.$ In particular,
$\bar{u}_{j+2}$ will be a quasi-regulous function and, by Lemma \ref{symmetr},
$\bar{u}_{j+2}, \bar{v}_{j+2}$ will be of class $\mathcal{C}^{k+1},$ which will complete the
proof.

From the input data we have semialgebraic functions $v_{j+1}, u_{j+2}, v_{j+2},$ {from
$\mathbf{R}^n\setminus\Sigma$ to $\mathbf{R},$} satisfying $v_{j+1}^2+u_{j+2}^2=v_{j+2}^2.$
Let us denote
$$2\alpha=\sqrt{v_{j+2}-v_{j+1}}(\sqrt{v_{j+2}-v_{j+1}}+\sqrt{v_{j+2}+v_{j+1}})+v_{j+1}.$$
This is well defined {on $\mathbf{R}^n\setminus\Sigma$} as, by the definition of $v_{j+2},$
we have $v_{j+2}\geq v_{j+1}\geq 0$. Observe that $2\alpha=v_{j+2}+|u_{j+2}|$ so
$|u_{j+2}|=\alpha-\frac{v_{j+1}^2}{4\alpha}$ and $v_{j+2}=\alpha+\frac{v_{j+1}^2}{4\alpha}.$

{From the equation $v_{j+1}^2+u_{j+2}^2=v_{j+2}^2$ and from the fact that $v_{j+2}\geq
v_{j+1}\geq 0$ it follows that $(\sqrt{v_{j+2}-v_{j+1}})^{-1}(0)=u_{j+2}^{-1}(0).$ Let
$f:\mathbf{R}^n\rightarrow\mathbf{R}$ be a continuous nonnegative semialgebraic function with
the following properties: $f|_M=(\sqrt{v_{j+2}-v_{j+1}})|_M,$
$f^{-1}(0)=(\sqrt{v_{j+2}-v_{j+1}})^{-1}(0)\cup\Sigma.$} Such a function can be obtained in
the same way as in the first step.

{Since $f^{-1}(0)=u_{j+2}^{-1}(0)\cup\Sigma$ is a compact $\mathcal{C}^{\infty}$ submanifold
of $\mathbf{R}^n$ of pure dimension $n-1$, there is an open ball
$G\subset\subset\mathbf{R}^n$ such that $f^{-1}(0)\cup M\subset G.$ Then, $f|_{\overline{G}}$
and $g|_{\overline{G}}=1$ satisfy the hypotheses of Lemma \ref{level}.}

By Lemma \ref{level}, we obtain polynomial approximation $\bar{\gamma}$ of {$f$ on
$\overline{G}$} of the form $\bar{\gamma}=\tilde{f}\cdot P^{2l}$ with any given precision
$\varepsilon$ (we assume $\varepsilon<1$), where $P, \tilde{f}\in\mathbf{R}[x_1,\ldots,x_n],$
$\tilde{f}>0$ on {$\mathbf{R}^n$} and $l=4^{m-j}(k+1)^{2(m-j+1)}.$ In particular, we have
{$\bar{\gamma}\in\mathcal{C}^{k+1}_{{4^{m-j}(k+1)^{2(m-j+1)}}}(\mathbf{R}^n)$}. Next
approximate, using the Stone-Weierstrass approximation theorem, the function
$(\sqrt{v_{j+2}-v_{j+1}}+\sqrt{v_{j+2}+v_{j+1}})$ {on $M$} by a polynomial $\bar{A}$ positive
on {$\mathbf{R}^n$} and define $2\bar{\alpha}:=\bar{\gamma}\bar{A}+\bar{v}_{j+1}$ {on
$\mathbf{R}^n$.} Finally define
$|\bar{u}_{j+2}|:=\bar{\alpha}-\frac{\bar{v}_{j+1}^2}{4\bar{\alpha}}$ and
$\bar{v}_{j+2}:=\bar{\alpha}+\frac{\bar{v}_{j+1}^2}{4\bar{\alpha}}.$

Using Lemma \ref{descent} we obtain
{$\frac{\bar{v}_{j+1}^2}{4\bar{\alpha}}\in\mathcal{C}^{k+1}_{2\cdot
4^{m-j-1}(k+1)^{2(m-j+1)-1}}(\mathbf{R}^n)$ and $\bar{v}_{j+2}\in\mathcal{C}^{k+1}_{
4^{m-j-1}(k+1)^{2(m-j)}}(\mathbf{R}^n).$} We also have
{$|\bar{u}_{j+2}|\in\mathcal{C}^{k+1}_{ 4^{m-j-1}(k+1)^{2(m-j)}}(\mathbf{R}^n).$} Indeed,
$|\bar{u}_{j+2}|=\bar{\alpha}-\frac{\bar{v}_{j+1}^2}{4\bar{\alpha}}=$
$$
\frac{4\bar{\alpha}^2-\bar{v}_{j+1}^2}{4\bar{\alpha}}=
\frac{\bar{\gamma}\bar{A}(\bar{\gamma}\bar{A}+2\bar{v}_{j+1})}{2(\bar{\gamma}\bar{A}+\bar{v}_{j+1})}=
\frac{\bar{\gamma}\bar{A}}{2}\cdot
(1+\frac{\bar{v}_{j+1}}{\bar{\gamma}\bar{A}+\bar{v}_{j+1}})=\frac{\bar{\gamma}\bar{A}}{2}+\frac{1}{2}\cdot\frac{\bar{\gamma}\bar{A}\bar{v}_{j+1}}{\bar{\gamma}\bar{A}+\bar{v}_{j+1}}.$$
Now the assertion follows by Lemma \ref{descent}.

It is clear that $\bar{v}_{j+2}$ and $|\bar{u}_{j+2}|$ defined above are nonnegative regulous
functions {on $\mathbf{R}^n$} satisfying $\bar{v}_{j+1}^2+|\bar{u}_{j+2}|^2=\bar{v}_{j+2}^2.$
{Note that $|\bar{u}_{j+2}|^{-1}(0)=(\bar{\gamma})^{-1}(0)$ is an algebraic subset of
$\mathbf{R}^n$.} It is also not difficult to observe that if {$\bar{\gamma}|_{M},
\bar{A}|_{M}$ and $\bar{v}_{j+1}|_{M}$ approximate $\sqrt{v_{j+2}-v_{j+1}}|_{M}=f|_{M},$
$(\sqrt{v_{j+2}-v_{j+1}}+\sqrt{v_{j+2}+v_{j+1}})|_{M}$ and $v_{j+1}|_{M},$} respectively,
with precision high enough, then $\bar{v}_{j+2},$ $|\bar{u}_{j+2}|$ are good approximations
of $v_{j+2},$ $|u_{j+2}|$ {on ${M}$}.

To complete the step number $j+1$ it remains to recover $\bar{u}_{j+2}$ approximating
$u_{j+2}$ {on $M$} when $|\bar{u}_{j+2}|$ approximating $|u_{j+2}|$ is given. To do this, it
is sufficient to define the sign of $\bar{u}_{j+2}$ on every connected component of
{$\mathbf{R}^n\setminus |\bar{u}_{j+2}|^{-1}(0)$} and to check that the obtained function
satisfies the requirements.

From the equation $|\bar{u}_{j+2}|=\frac{\bar{\gamma}\bar{A}}{2}\cdot
(1+\frac{\bar{v}_{j+1}}{\bar{\gamma}\bar{A}+\bar{v}_{j+1}})$ above we know that the zeroes of
$|\bar{u}_{j+2}|$ are precisely the zeroes of $\bar{\gamma}$ (where $\bar{\gamma}$ is the
polynomial approximating {$f$ on $\overline{G}$} obtained by Lemma \ref{level}). Let
$A_1,\ldots,A_s$ denote all connected components of {$G\setminus f^{-1}(0).$} By
Lemma~\ref{level} {(recall that $g:=1$)}, there is a tubular neighborhood $T\subset\subset G$
of {$f^{-1}(0)$} and there are connected components $B_1,\ldots,B_s$ of
{$G\setminus(\bar{\gamma})^{-1}(0)$} with $A_i\subset B_i\cup\overline{T}$ and such that
$B_i\neq B_j$ if $i\neq j.$ (Here recall that we assume $\varepsilon<1$ so the polynomial
approximating $g,$ obtained by applying Lemma~\ref{level}, does not vanish at any point of
{$\overline{G}$}.)

Recall that {$u_{j+2}^{-1}(0)\cup\Sigma=f^{-1}(0)$} and set
$\mathrm{sgn}(\bar{u}_{j+2}|_{B_i})=\mathrm{sgn}(u_{j+2}|_{A_i}),$ for $i=1,\ldots,s.$ On
every component of {$G\setminus(\bar{\gamma})^{-1}(0)$} different from $B_1,\ldots, B_s$ set
the sign of $\bar{u}_{j+2}$ arbitrarily and set $\bar{u}_{j+2}|_{(\bar{\gamma})^{-1}(0)}=0$.
{Since $(\bar{\gamma})^{-1}(0)\cap\partial G=\emptyset,$ there is precisely one connected
component $C$ of $G\setminus(\bar{\gamma})^{-1}(0)$ with $\overline{C}\cap\partial
G\neq\emptyset.$ Therefore on every connected component $E$ of $(\mathbf{R}^n\setminus
G)\setminus(\bar{\gamma})^{-1}(0)$ we can define
$\mathrm{sgn}(\bar{u}_{j+2}|_{E})=\mathrm{sgn}(\bar{u}_{j+2}|_C$).} Clearly, $\bar{u}_{j+2}$
is a quasi-regulous function and, by Lemma \ref{symmetr},
$\bar{u}_{j+2}\in\mathcal{C}^{k+1}(\mathbf{R}^n).$

Let us check that $\bar{u}_{j+2}$ approximates $u_{j+2}$ on {$B_i\cap{M},$} for
$i=1,\ldots,s.$ It is clear that $\bar{u}_{j+2}$ approximates $u_{j+2}$ on {$A_i\cap B_i\cap
M.$} This is because $u_{j+2}$ and $\bar{u}_{j+2}$ have the same sign at every point of
${A_i\cap B_i}$ and $|\bar{u}_{j+2}|$ approximates $|u_{j+2}|$ {on $M.$} On the other hand,
$B_i\setminus A_i\subset\overline{T}.$ If this was not true, then we had $((B_i\setminus
A_i)\setminus\overline{T})\cap A_j\neq\emptyset$ for some $j\neq i$ because
{$\bigcup_{i=1}^sA_i\cup\overline{T}=G.$} Then $((B_i\setminus A_i)\setminus\overline{T})\cap
(B_j\cup\overline{T})\neq\emptyset$ so $B_i\cap B_j\neq\emptyset$ for some $j\neq i,$ a
contradiction. Since $B_i\setminus A_i\subset\overline{T},$ the function $\bar{u}_{j+2}$
approximates $u_{j+2}$ on {$(B_i\setminus A_i)\cap{M}$} because
{$|u_{j+2}|_{M\cap\overline{T}}|=(f\cdot\sqrt{v_{j+2}+v_{j+1}})|_{M\cap\overline{T}},$
$|\bar{u}_{j+2}|_{M\cap\overline{T}}|$ approximates $|u_{j+2}|_{M\cap\overline{T}}|$} and, by
Lemma \ref{level}, {$(f\cdot\sqrt{v_{j+2}+v_{j+1}})|_{M\cap\overline{T}}$ may be assumed as
close to zero as we wish. Hence, $u_{j+2}|_{M\cap\overline{T}}$ and
$\bar{u}_{j+2}|_{M\cap\overline{T}}$} are close to zero and therefore close to each other
(even if $u_{j+2}$ and $\bar{u}_{j+2}$ have different signs at some points of
{$M\cap\overline{T}$}).

Let us check that $\bar{u}_{j+2}$ approximates $u_{j+2}$ on
{$M\setminus\bigcup_{i=1}^s{B_i}.$} Observe that {$G\setminus\bigcup_{i=1}^s
{B_i}\subset\overline{T}.$} This is because {$\bigcup_{i=1}^s{A_i}\cup\overline{T}={G}$} and,
for every $i,$ we have ${A_i}\subset{B_i}\cup \overline{T}.$ Consequently,
{$M\setminus\bigcup_{i=1}^s{B_i}\subset\overline{T}$} and the assertion follows as above.
Thus we have proved that $\bar{u}_{j+2}$ approximates $u_{j+2}$ on {$M.$}

{Finally, recall that, by construction, $\bar{u}_j^{-1}(0)$ is an algebraic subset of
$\mathbf{R}^n,$ for every $j=1,\ldots,m+1.$ Therefore,
$\bar{v}_{m+1}^{-1}(0)=\bar{u}_1^{-1}(0)\cap\ldots\cap\bar{u}_{m+1}^{-1}(0)$ is an algebraic
subset of $\mathbf{R}^n$. If the approximation is close enough, then
$\bar{v}_{m+1}^{-1}(0)\cap M=\emptyset$ because $v_{m+1}=1.$ In particular, every
$w_j=\frac{\bar{u}_j}{\bar{v}_{m+1}}|_{M}$ is a $\mathcal{C}^k$ quasi-regulous function
approximating $u_j|_M$ and $w_1^2+\ldots+w_{m+1}^2=1$.} \qed\vspace*{3mm}\\
{\small\textit{Acknowledgements.} I express my gratitude to Professor Wojciech
Kucharz, Professor Krzysztof Nowak and Professor Wies\l aw Paw\l ucki for comments, remarks
and very helpful discussions.}}

\end{document}